\documentclass[12pt]{amsart}

\usepackage{amsfonts,amsmath,amssymb,enumerate,
hyperref,latexsym,xcolor,amsthm,graphicx,multirow, threeparttable}
\def\a{\alpha}
\def\b{\beta}
\def\l{\langle} \def\r{\rangle}
\def\div{\,\,\big|\,\,}
\newcommand\Ga{\Gamma}
\def\O{{\bf O}}

\newcommand\ZZ{\mathbb{Z}} 
 
\newcommand\Aut{\mathrm{Aut}}  
  
\newcommand\soc{\mathrm{soc}}
\newcommand\C{\mathbf{C}}\newcommand\CC{\mathcal{C}}

\newcommand\GL{\mathrm{GL}}

\newcommand\GammaL{\mathrm{\Gamma L}}

 \newcommand\Q{\mathrm{Q}}

\newtheorem{theorem}{Theorem}[section]

\newtheorem{corollary}[theorem]{Corollary}

\newtheorem{lemma}[theorem]{Lemma}

\theoremstyle{definition}

\begin{document}
\title[2-arc-transitive graphs]{On basic $2$-arc-transitive graphs}
\thanks{Supported by the National Natural Science Foundation of China
(11971248, 11731002), and the Fundamental Research Funds for the Central
Universities.}

\author[Z.P. Lu]{Zai Ping Lu}
\address{Z.P. Lu\\ Center for Combinatorics, LPMC\\ Nankai University\\ Tianjin 300071\\ P. R. China}
\email{lu@nankai.edu.cn}

\author[R.Y. Song]{Ruo Yu Song}
\address{R.Y. Song\\ Center for Combinatorics, LPMC\\ Nankai University\\ Tianjin 300071\\ P. R. China}
\email{851878536@qq.com}

\begin{abstract}
A connected graph $\Ga=(V,E)$ of valency at least $3$ is called a basic $2$-arc-transitive graph if
its full automorphism group has a subgroup $G$ with the following
properties: (i) $G$ acts transitively on the set of $2$-arcs of $\Ga$, and (ii) every minimal normal subgroup of $G$ has at most two orbits on $V$.
 Based on Praeger's theorems on  $2$-arc-transitive graphs, this paper presents  a  further understanding on the automorphism group of a basic $2$-arc-transitive graph.

\vskip 5pt

\end{abstract}
\maketitle

\baselineskip 15pt


\section{Introduction}
All graphs considered in this paper are assumed to be finite, simple
and undirected.

\vskip 5pt

Let $\Ga=(V,E)$ be a graph with vertex set $V$ and edge set $E$. Denote by $\Aut(\Ga)$ the  full automorphism group of the graph $\Ga$.
A subgroup $G$ of $\Aut(\Ga)$, written as $G\leqslant \Aut(\Ga)$, is called a  group of $\Ga$.
For a vertex $\a\in V$, let $G_\a=\{g\in G\mid \a^g=\a\}$ and $\Ga(\a)=\{\b\in V\mid \{\a,\b\}\in E\}$, called the stabilizer of $\a$ in $G$ and the neighborhood of $\a$ in $\Ga$, respectively. A group $G$ of $\Ga$ is call locally-primitive on $\Ga$ if
for each $\a\in V$ the stabilizer $G_\a$ acts primitively on $\Ga(\a)$,
that is, $\Ga(\a)$ has no nontrivial $G_\a$-invariant partition.
Recall that an arc of $\Ga$ is an ordered pair of adjacent vertices, and a $2$-arc
is a triple $(\alpha,\beta,\gamma)$  of  vertices with $\{\a, \b\},\,\{\b, \gamma\}\in E$ and
$\a\neq\gamma$.
A  group $G$ of $\Ga$ is said to be vertex-transitive, edge-transitive, arc-transitive or $2$-arc-transitive on $\Ga$ if $G$ acts transitively on the vertices, edges, arcs or $2$-arcs of $\Ga$, respectively. A graph is called vertex-transitive, edge-transitive, arc-transitive or $2$-arc-transitive if it has a vertex-transitive, edge-transitive, arc-transitive or $2$-arc-transitive group, respectively.

\vskip 5pt

A connected regular graph $\Ga=(V,E)$ of valency at least $3$ is called a basic $2$-arc-transitive graph if it has a $2$-arc-transitive  group $G$ such that every minimal normal subgroup of $G$ has at most two orbits on $V$.
Praeger \cite{Praeger-ON,Prag-o'Nan-bi} observed that a connected $2$-arc-transitive graph of valency at least $3$ is a normal cover of some basic $2$-arc-transitive graph.
Based on the O'Nan-Scott theorem for quasiprimitive permutation groups established in \cite{Praeger-ON},
Praeger \cite{Praeger-ON,Prag-o'Nan-bi} characterized the group-theoretic structures
for  basic $2$-arc-transitive graphs. She  proved that, except for one case about bipartite graphs, basic $2$-arc-transitive graphs are associated with  quasiprimitive groups
of type I, II, IIIb(i) or III(c) described as in \cite[Section 2]{Praeger-ON}, which is named  HA, AS, PA or TW  in \cite{Prag-o'Nan-survey}, respectively.

\vskip 5pt

Praeger's framework for $2$-arc-transitive graphs stimulated a wide interest in  classification or characterization of basic $2$-arc-transitive graphs.
For example, a construction of the graphs associated with quasiprimitive permutation groups of type TW is given in \cite{Baddeley}, the graphs associated with Suzuki simple groups, Ree simple groups and 2-dimensional projective
linear groups are classified  in \cite{FP1,FP2,HNP} respectively, the graphs of order a prime power are classified in \cite{prime-power}.
Besides, Li \cite{Li1} proved that all basic $2$-arc-transitive graphs of odd order can be constructed from almost simple groups, which inspires the ongoing project to classify basic $2$-arc-transitive graphs of odd order, see \cite{odd-An} for some progress in this topic.

\vskip 5pt

In this paper, we have a further understanding on the automorphism groups of basic $2$-arc-transitive graphs, which may be helpful to study the Praeger's problem proposed in \cite{Prag-o'Nan-bi}: {\em Classify all finite basic $2$-arc-transitive graphs}.
Assume that $\Ga=(V,E)$ is a basic $2$-arc-transitive graph with respect to $G$. Fix an edge $\{\a,\b\}\in E$, and set $G^*=\l G_\a,G_\b\r$. It is well-known that $|G:G^*|\leqslant 2$, $G^*$ is edge-transitive on $\Ga$, and $\Ga$ is bipartite if and only if $|G:G^*|=2$,
refer to \cite[Exercise 3.8]{Weiss}.
If $\Ga$ is not bipartite, then $G$ is a quasiprimitive permutation  group on $V$ of type HA, AS, PA or TW,  refer to \cite[Theorem 2]{Praeger-ON} or \cite[Theorem 6.1]{Prag-o'Nan-survey}; in this case, $G$ has a unique minimal normal subgroup. Now assume that $\Ga$ is bipartite, that is, $|G:G^*|=2$. Praeger \cite{Prag-o'Nan-bi} proved that
either $\Ga$ is a complete bipartite graph, or $G^*$ is faithful on both parts of $\Ga$ and one of the following holds:
\begin{enumerate}[\rm (I)]
  \item   $G^*$ is  quasiprimitive on both parts of $\Ga$  with a same type HA, AS, PA or TW;
  \item   $G$ has a  normal subgroup $N$ which is a direct product of two intransitive minimal normal subgroups of $G^*$.
\end{enumerate}
In the paper, we prove that the normal subgroup $N$ in (II)  is
 the unique minimal normal subgroups of $G$, see Theorem \ref{basic-gph-1}.
Then we  formulate the following  result, which is finally proved in Section \ref{sect=main-pf}.
\begin{theorem}\label{basic-gph}
Assume that $\Ga=(V,E)$ is a basic $2$-arc-transitive graph with respect to a group $G$. Let $G^*=\l G_\a,G_\b\r$ and $N=\soc(G^*)$, where  $\{\a,\b\}\in E$.
Then either $\Ga$ is a complete bipartite graph, or the  following statements hold:
\begin{itemize}
\item[(1)]   $N$  is the unique minimal normal subgroup of $G$;
  \item[(2)]  either  $N$ is simple, or every simple direct factor of $N$ is semiregular on $V$;
  \item[(3)]   either $N$ is locally-primitive on  $\Ga$, or $N_\a$ is given as follows:
\begin{itemize}
        \item[(i)] $N_\a=\ZZ_p^k{:}(\ZZ_{m_1}{.}\ZZ_m)=(\ZZ_p^k\times \ZZ_{m_1}){.}\ZZ_m$ and $|\Ga(\a)|=p^k$, where $m_1\div m$, $m\div (p^d-1)$ for some divisor $d$ of $k$ with $d<k$; or
        \item[(ii)] $N_\a=\ZZ_3^4{:}(Q.\Q_8)=(\ZZ_3^4\times Q){.}\Q_8$ and $|\Ga(\a)|=3^4$, where $\Q_8$ is the quaternion group and $Q$ is isomorphic to a subgroup of $\Q_8$.
\end{itemize}
\end{itemize}
\end{theorem}

\vskip 5pt

It is well-known that the order of a finite nonabelian simple group is divisible by $4$ and two distinct odd primes. By (2) of Theorem \ref{basic-gph}, we have the following corollaries.

\begin{corollary}\label{cor-1}
Assume that $\Ga=(V,E)$ is a  basic $2$-arc-transitive graph
with respect to a group $G$, and   $G^*=\l G_\a,G_\b\r$ for an edge $\{\a,\b\}\in E$. Assume that one of $G^*$-orbits on $V$ has length
$p^aq^b$, where $a$ and $b$ are  positive integers, $p$ and $q$ are distinct primes. If  $\Ga$ is not a complete bipartite graph, then $G$ is almost simple.
\end{corollary}

\begin{corollary}\label{cor-2}
Assume that $\Ga=(V,E)$ is a  basic $2$-arc-transitive graph
with respect to a group $G$, and   $G^*=\l G_\a,G_\b\r$ for an edge $\{\a,\b\}\in E$.  Assume that one of $G^*$-orbits on $V$ has length
$n$ or $2n$, where $n$ is either an odd integer or a power of $2$. If $\Ga$ is not a complete bipartite graph then either $G$ is almost simple, or $n=|G^*:G_\a|=p^k$ and $\soc(G^*)\cong \ZZ_p^k$, where $p$ is a  prime and $k\geqslant 1$.
\end{corollary}

Another consequence of Theorem \ref{basic-gph} is stated as follows.

\begin{theorem}\label{bi-prim}
Let $\Ga=(V,E)$ be a connected graph, $G\leqslant \Aut(\Ga)$ and $G^*=\l G_\a,G_\b\r$ for an edge $\{\a,\b\}\in E$.
Assume that $G$ is $2$-arc-transitive on $\Ga$, and $G^*$ acts primitively on
each of $G^*$-orbits on $V$. Then one of the following holds:
\begin{itemize}
\item[(1)]   $\Ga$  is a complete bipartite graph;
\item[(2)]  $\soc(G)=\soc(G^*)$, and $\soc(G^*)$ is either simple or regular on each $G^*$-orbit;
  \item[(3)]  $\Ga$ is bipartite, $\soc(G)=\soc(G^*)\times M$ with $|M|=2$, and $\soc(G^*)$ is either simple or regular on each $G^*$-orbit.
\end{itemize}
\end{theorem}

\vskip 20pt

\section{Some observations on $2$-transitive permutation  groups}\label{sect=2tg}
This section gives some simple results about $2$-transitive permutation groups, which serves to analyze the structures of vertex-stabilizers of $2$-arc-transitive graphs.

\vskip 5pt

Let $X$ be a transitive permutation group on a finite set $\Omega$.
Recall that the socle $\soc(X)$ is generated by all minimal normal subgroups of $X$. It is easily shown that  $\soc(X)$ is a characteristic subgroup of $X$.
Assume that $X$ is a $2$-transitive permutation group on   $\Omega$.
Then   $\soc(X)$ is either elementary abelian and
regular $\Omega$, or simple and primitive on $\Omega$, refer to \cite[page 101, Theorem 4.3]{Cameron} and \cite[page 107, Theorem 4.1B]{Dixon}. In particular, $X$ is either affine or almost simple. Inspecting  the list of finite $2$-transitive permutation groups
(refer to \cite[pages 195-197, Tables 7.3 and 7.4]{Cameron}), we have the following basic fact, see also \cite[Corollary 2.5]{LSS}.

\vskip 5pt

\begin{lemma}\label{stab-2tg}
Let $X$ be a $2$-transitive permutation group on a finite set $\Omega$, and $\a\in \Omega$. Assume that $K$ is an insolvable normal subgroup of $X_\a$. Then $K$ has a unique insolvable composition factor say $S$, and $S$ is isomorphic to a composition factor of $X$
if and only if $X$ is affine.
\end{lemma}

Recall that a transitive permutation group $X$  on $\Omega$ is  a Frobenius group
if $X$ is not regular on $\Omega$ and, for $\a\in \Omega$, the point-stabilizer $X_{\a}$, called a Frobenius complement of $X$, is semiregular on  $\Omega$.

 \vskip 5pt

By Frobenius' Theorem (refer to \cite[pages 190-191, (35.23) and (35.24)]{Aschbacher}), for a Frobenius group $X$ on $\Omega$, the identity and elements
without fixed-point form a normal regular subgroup of $X$, which is called  the Frobenius kernel of $X$.

\begin{lemma}\label{Frobenius}
Let $X=KH$ be an imprimitive Frobenius group on $\Omega$ with Frobenius kernel $K\cong \ZZ_p^k$ and a Frobenius complement $H$, where $p$ is a prime  and   $k\geqslant 2$. Then $H$ is isomorphic to an irreducible subgroup of the general linear group $\GL_l(p)$, and $|H|$ is a divisor of $p^d-1$, where $2l\leqslant k$ and
$d$ is a common divisor of $k$ and $l$.
 \end{lemma}
\proof
Note that $H$ acts faithfully and semiregularly on $K\setminus \{1\}$ by conjugation, see \cite[page 191, (35.25)]{Aschbacher}. Then $|H|$ is a divisor of $p^k-1$.
Assume that $X$ is imprimitive on $\Omega$. Then $K$ is not a minimal normal subgroup of $X$. By Maschke's Theorem (refer to \cite[page 40, (12.9)]{Aschbacher}),
$K$ is a direct product of two $H$-invariant proper subgroups.
Thus we may choose a minimal $H$-invariant   subgroup $L$ of $K$ with $|L|^2\leqslant |K|$. It is easily shown that $LH$ is a primitive Frobenius group (on an $L$-orbit), which has  Frobenius kernel $L$.
Set $|L|=p^l$. Then  $|H|$ is a divisor of $p^l-1$,
$2l\leqslant k$, and $H$ is isomorphic to an irreducible subgroup of $\GL_l(p)$.

\vskip 5pt

Choose a minimal positive integer $d$ such that $|H|$ is a divisor of $p^d-1$. Then $d\leqslant l$.
Set $k=xd+y$ for integers $x\geqslant 1$ and $0\leqslant y<d$. Then $q^k-1=q^y(q^{xd}-1)+(q^y-1)$, and thus
$|H|$ is a   divisor of $q^y-1$. By the choice of $d$, we have $y=0$, and so $d$
is a divisor of $k$. Similarly, $d$ is a divisor of $l$. Then the lemma follows.
\qed

\begin{lemma}\label{imp-normal-sub}
Let $X$ be a $2$-transitive permutation group on a finite set $\Omega$.
Assume that $1\ne  N\unlhd X$. Then $\soc(N)=\soc(X)$, and either
$N$ is primitive on $\Omega$ or
one of the following holds:
 \begin{enumerate}[\rm (1)]\itemindent=1em
        \item  $N=\ZZ_p^k{:}\ZZ_m$ and $|\Omega|=p^k$, where $p$ is a prime, $k\geqslant 2$, $m\div (p^d-1)$ for some divisor $d$ of $k$ with $d<k$;
        \item  $N=\ZZ_3^4{:}\Q_8$ and $|\Omega|=3^4$.
\end{enumerate}
\end{lemma}
\proof
Noting that $\soc(N)$ is characteristic in $N$, it follows that $\soc(N)$ is a normal subgroup of $X$.
Since $\soc(X)$ is the unique minimal normal subgroup of $X$,
we have $\soc(X)\leqslant \soc(N)$.
Suppose that $\soc(X)\ne \soc(N)$. Then
$\soc(N)$ has a simple direct factor $T$ with $T\cap \soc(X)=1$. This implies that
$T$ lies in the centralizer $\C_{X}(\soc(X))$ of $\soc(X)$; however, by \cite[page 109, Theorem 4.2A]{Dixon}, $\C_{X}(\soc(X))=1$ or $\soc(X)$,
a contradiction. Then $\soc(N)=\soc(X)$.

\vskip 5pt

Next we  assume that $N$ is imprimitive
on $\Omega$, and show one of (1) and (2) holds.
By  \cite[pages 215-217, Theorems 7.2C and 7.2E]{Dixon},
$\soc(N)=\soc(X)\cong \ZZ_p^k$ for a prime $p$ and integer $k\geqslant 2$ with $|\Omega|=p^k$, and either $N=\soc(X)$   or $N$ is a
Frobenius group with Frobenius kernel $\soc(X)$.
In particular, by Lemma \ref{Frobenius}, we write
$N=KH$, where $K\cong \ZZ_p^k$ and $|H|$ is a divisor of
$p^d-1$ for a divisor $d$ of $k$ with $d<k$.
Note that $X$ is an affine $2$-transitive permutation group.
Inspecting the affine $2$-transitive permutation groups listed in \cite[page 197, Table 7.4]{Cameron}, we conclude that either
\begin{itemize}
  \item[(i)] either $H$ is cyclic, or $H\leqslant \GammaL_1(p^k)$; or
  \item [(ii)] $p^k=3^4$, yielding $d\in \{1,2\}$, and so $|H|$ is a divisor of $8$.
\end{itemize}

Suppose that (i) holds. Note that $H$ acts faithfully and semiregularly on $K\setminus \{1\}$ by conjugation. If $H\leqslant \GammaL_1(p^k)$ then $H\leqslant\GL_1(p^k)$, otherwise,
$H$ should contain an element which centralizes some element in $K\setminus \{1\}$, a contradiction. Thus $H$ is cyclic, and so $H\lesssim\ZZ_{p^d-1}$. Then part (1) of this lemma follows.

\vskip 5pt

Suppose that (ii) holds. If $H$ is cyclic then $N$ is described as in part
(1) of this lemma. Assume that $H$ is not cyclic. By Lemma \ref{Frobenius}, $H$ is isomorphic to an irreducible subgroup of  $\GL_2(3)$.
Confirmed by GAP \cite{GAP}, we have  $H\cong \Q_8$, and thus part (2) of this lemma follows. This completes the proof.
\qed

\begin{lemma}\label{affine-2trans}
Let $X$ be an affine $2$-transitive permutation group, and $\soc(X)=K_1\times
\cdots \times K_l$, where $1<K_i<\soc(X)$ for $1\leqslant i\leqslant l$.
Then there are $x\in X$ and $i$ such that $K_i^x\not\in\{ K_i\mid 1\leqslant i\leqslant l\}$.
\end{lemma}
\proof
Clearly, $\cup_{i}(K_i\setminus\{1\})\ne \soc(X)\setminus\{1\}$.
Let $H$ be a point-stabilizer of $X$. Then $H$ acts transitively on $\soc(X)\setminus\{1\}$ by conjugation. Thus $H$ does not fix  $\cup_{i}(K_i\setminus\{1\})$ set-wise by conjugation, and the lemma follows.
\qed

\vskip 20pt

\section{The uniqueness of minimal normal subgroup}\label{sect=local}

In this section, we assume that $\Ga=(V,E)$ is a connected  regular graph, and $G\leqslant \Aut(\Ga)$.
Denote by $G_\a^{\Ga(\a)}$ the permutation group induced by $G_\a$ on
$\Ga(\a)$. Let $G_\a^{[1]}$ be the kernel of $G_\a$ acting on $\Ga(\a)$. Then
\[G_\a^{\Ga(\a)}\cong G_\a/G_{\a}^{[1]}.\]

\vskip 5pt

Let $\b\in \Ga(\a)$, and set $G_{\a\b}^{[1]}=G_\a^{[1]}\cap G_\b^{[1]}$.
Then $G_{\a\b}^{[1]}$ is the kernel of the arc-stabilizer $G_{\a\b}$ acting
on $\Ga(\a)\cup\Ga(\b)$. Noting that $G_\a^{[1]}\unlhd G_{\a\b}$, we have
\[
G_\a^{[1]}/G_{\a\b}^{[1]}\cong (G_\a^{[1]})^{\Ga(\b)}\unlhd G_{\a\b}^{\Ga(\b)}=(G_{\b}^{\Ga(\b)})_\a.\]

\vskip 5pt

Assume that  $G$ is arc-transitive on $\Ga$, and $N$ is an arbitrary  normal subgroup of $G$. Then \[N_\a\unlhd G_\a,\,N_\a^{[1]}\unlhd G_\a^{[1]},\, N_{\a\b}\unlhd G_{\a\b},\, N_{\a\b}^{[1]}\unlhd G_{\a\b}^{[1]}.\]
Taking $x\in G$ with $(\a,\b)^x=(\b,\a)$, we have  \[N_\b=N_\a^x,\, N_{\a\b}^x=N_{\a\b},\,\Ga(\b)=\Ga(\a)^x.\] It follows that $N_{\a\b}^{\Ga(\b)}\cong N_{\a\b}^{\Ga(\a)}$.
Since  $N_{\a\b}\unlhd G_{\a\b}$, we have $N_{\a\b}^{\Ga(\a)}\unlhd G_{\a\b}^{\Ga(\a)}$, and so
\begin{equation}\label{ext-2}
N_\a^{[1]}/N_{\a\b}^{[1]}\cong (N_\a^{[1]})^{\Ga(\b)}\unlhd N_{\a\b}^{\Ga(\b)}\cong N_{\a\b}^{\Ga(\a)}=(N_{\a}^{\Ga(\a)})_\b\unlhd (G_{\a}^{\Ga(\a)})_\b.
\end{equation}
In particular, $N_\a^{[1]}/N_{\a\b}^{[1]}$ is isomorphic to a normal subgroup
of $(N_{\a}^{\Ga(\a)})_\b$.

\vskip 5pt

Assume that $G$ is $2$-arc-transitive on $\Ga$. Then   $G_{\a\b}^{[1]}$ has order a prime power, see \cite[Corollary 2.3]{Gardiner}. In particular, $G_{\a\b}^{[1]}$ is solvable.
Then $(G_\a^{[1]})^{\Ga(\b)}$ is solvable if and only if $G_\a^{[1]}$ is solvable.
Note that $G_\a^{\Ga(\a)}$ is a $2$-transitive group on $\Ga(\a)$.
By Lemma \ref{stab-2tg},
both $G_{\a}^{\Ga(\a)}$ and $(G_{\a}^{\Ga(\a)})_\b$ have at most one insolvable composition factor.
Then we have the following fact.

\begin{lemma}\label{factors}
Assume that $G$ is $2$-arc-transitive on $\Ga$, $N\unlhd G$  and $N_\a$ is insolvable, where $\a\in V$. Then
$N_\a^{\Ga(\a)}$ has a unique insoluble composition factor, and $N_\a^{[1]}$ has at most one insoluble composition factor. If $N_\a^{[1]}$ and $N_\a^{\Ga(\a)}$ have isomorphic insoluble composition factors then $G_\a^{\Ga(\a)}$ is an affine $2$-transitive permutation group.
\end{lemma}

\begin{lemma}\label{simple-stab}
Let $N\unlhd G$ and $\a\in V$.
Assume that $G$ is $2$-arc-transitive on $\Ga$, and that $N_\a$ has a normal subgroup $K\cong T^k$ for an integer $k\ge 1$ and a nonabelian simple group $T$.
Then $k=1$.
\end{lemma}
\proof
Note that every normal subgroup of $K$ is isomorphic to $T^l$ for some $l\leqslant k$, where $T^0=1$. Set $K\cap G_\a^{[1]}\cong T^l$. Then
\[K^{\Ga(\a)}\cong KG_\a^{[1]}/G_\a^{[1]}\cong K/(K\cap G_\a^{[1]})\cong T^{k-l}.\] Since $K^{\Ga(\a)}\unlhd N_\a^{\Ga(\a)}\unlhd  G_\a^{\Ga(\a)}$,
by   Lemma \ref{factors},
we conclude that $l,k-l\in \{0,1\}$.
If $G_\a^{\Ga(\a)}$ is of affine type, then $k-l=0$, and so $k=l=1$.
If $G_\a^{\Ga(\a)}$ is almost simple, then either $k=l=1$ or $k-l=1$ and $l=0$, and so $k=1$. This completes the proof.
\qed

\vskip 10pt

\begin{theorem}\label{basic-gph-1}
Assume that $\Ga=(V,E)$ is a basic $2$-arc-transitive graph with respect to  $G$, and $G^*=\l G_\a,G_\b\r$ for $\{\a,\b\}\in E$. Then either $\Ga$ is a complete bipartite graph, or $\soc(G^*)$ is the unique minimal normal subgroup of $G$  and one of the following holds:
\begin{itemize}
 \item[(1)]   $\soc(G^*)$ is semiregular on $V$;
  \item[(2)]   $\soc(G^*)$ is a nonabelian simple group;
  \item[(3)]   $G^*$ is a quasiprimitive permutation group of type {\rm PA} on each  $G^*$-orbit   on $V$;
 \item[(4)] $\Ga$ is a bipartite graph, $G^*$ is faithful on each part of $\Ga$,
 $\soc(G^*)=M_1\times M_2$ for  minimal normal subgroups
$M_1$ and $M_2$ of $G^*$, and both $M_1$ and $M_2$ are
semiregular and intransitive on each part of $\Ga$.
\end{itemize}
\end{theorem}
\proof
If $\Ga$ is not bipartite  then $G=G^*$ and, by  \cite[Theorem 2]{Praeger-ON}, $G$ has a unique minimal normal subgroup, and  one of parts (1)-(3) follows.
Thus we assume that $\Ga$ is a bipartite graph with two parts $U$ and $W$.
In particular,  $|G:G^*|=2$. By \cite[Theorem 2.1]{Prag-o'Nan-bi}, either $\Ga$ is a complete bipartite graph, or $G^*$ is faithful  on  each of $U$ and $W$. In the following, we assume  that the latter case occurs.

\vskip 5pt

Let $K$ be an arbitrary minimal normal subgroup of $G$.
Suppose that  $K\not\leqslant G^*$. Then $K\cap G^*=1$  and $G=G^*K$, yielding $|K|=2$.
Since $K$ has at most two orbits on $V$, we have $|V|\leqslant 4$, which is impossible as $\Ga$ is   bipartite and of valency at least $3$.
Therefore, $K\leqslant G^*$. Let $K_1$ be a  minimal normal subgroup of $G^*$
with $K_1\leqslant K$, and let $x\in G\setminus G^*$. Then $K_1^x$ is also a  minimal normal subgroup of $G^*$.
Noting that $x^2\in G^*$, we have
$(K_1^x)^x=K_1^{x^2}=K_1$. It follows that $K_1K_1^x$ is normal in $G$.
Since $K_1^x\leqslant K^x=K$, we have $K=K_1K_1^x\leqslant \soc(G^*)$.
It follows that $\soc(G)\leqslant \soc(G^*)$.

\vskip 5pt

{\bf Case 1}.
Assume that $G^*$ is quasiprimitive
on $U$ and $W$. Then, by\cite[Theorem 2.3]{Prag-o'Nan-bi}, $\soc(G^*)$ is the unique minimal normal subgroup of $G^*$, and one of  parts (1)-(3) of Theorem \ref{basic-gph-1}  occurs. In this case, since $\soc(G)\leqslant \soc(G^*)$, we know that $\soc(G^*)$ is the unique  minimal normal subgroup of $G$.

\vskip 5pt

{\bf Case 2}.
Assume that $G^*$ is not quasiprimitive
on  one of $U$ and $W$, say $U$. Then $G^*$ has a minimal normal subgroup  $M$
which is intransitive on $U$. Let $x\in G\setminus G^*$.
Then $M^x$ is a minimal normal subgroup of $G^*$, and $M^x$ is intransitive on $W$.
Note that $MM^x$ is normal in $G$. Then $MM^x$ is transitive on both $U$ and $W$. It follows that $M\ne M^x$, and so $M\cap M^x=1$. Then $MM^x=M\times M^x$.
If $M$ is transitive on $W$ then $M^x$ is semiregular on $W$ by \cite[Theorem 4.2A]{Dixon}, and thus both $M$ and $M^x$ are regular on $W$, a contradiction.
Therefore, $M$ is intransitive on $W$. Similarly, $M^x$ is intransitive on $U$.
It follows from \cite[Lemma 5.1]{GLP} that both $M$ and $M^x$ are semiregular on $U$ and $W$.

\vskip 5pt

Set $N=MM^x$, and write $M=T_1\times\cdots T_k$, where $T_i$ are isomorphic simple groups. Then \[N=T_1\times\cdots T_k\times T_1^x\times\cdots T_k^x.\] Let $L$ be a minimal normal subgroup of $G$ with $L\leqslant N$.
Assume that $M\not\leqslant L$. Then $M\cap L=1$ as $M$ is a  minimal normal subgroup of $G^*$, and so $M^x\cap L=(M\cap L)^x=1$. Thus both $M$ and $M^x$ centralize $L$. Considering the action of $G^*$ on $U$ or $W$, by \cite[Theorem 4.2A]{Dixon}, $L$ is nonabelian. This forces that
every $T_i$ is a nonabelian simple group. It follows that $L$ contains  $T_i$ or $T_i^x$ for some $i$. Then  $M\cap L\ne 1$ or $M^x\cap L\ne 1$, a contradiction.
Then $M\leqslant L$, yielding $N=MM^x\leqslant L$, and so $N=L$.
Therefore, $N$ is a minimal normal subgroup of $G$. Clearly, $N\leqslant \soc(G^*)$.

\vskip 5pt

Suppose that $N\ne \soc(G^*)$.
Then $G^*$ has a minimal normal subgroup $M_1$ with
$M_1\cap N=1$. This implies that $\C_{G^*}(N)\not\leqslant N$.
Noting that $\C_{G^*}(N)$ is a normal subgroup of $G$, it follows that $\C_{G^*}(N)$ acts transitively
on both $U$ and $W$.
 Considering the action of $G^*$ on $U$, it follows from \cite[Theorem 4.2A]{Dixon} that
$N$ is not abelian,  $\C_{G^*}(N)\cong N$, and both
 $\C_{G^*}(N)$ and $N$ are regular on $U$.
Let $\a\in U$ and $X=\C_{G^*}(N)N$. Then $X$ is normal in $G$, and
$X=\C_{G^*}(N)X_\a$. We have \[X_\a\cong \C_{G^*}(N)N/\C_{G^*}(N)\cong N=T_1\times\cdots T_k\times T_1^x\times\cdots \times T_k^x.\] Then $2k=1$ by Lemma \ref{simple-stab},  a contradiction.
Therefore, $N=\soc(G^*)$, and our result  follows.
\qed

\vskip 20pt

\section{Semiregular direct factors}\label{sect=main-pf}

Let $\Ga=(V,E)$ be a connected  regular graph, and $G\leqslant \Aut(\Ga)$.

\vskip 5pt

Assume that $G$ is a $2$-arc-transitive group of $\Ga$. Then
$G_\a^{\Ga(\a)}$ is a $2$-transitive permutation group on $\Ga(\a)$, where $\a\in V$. Let $N\unlhd G$ with
$N_\a\ne 1$. Then  it is easily shown that $N_\a$ acts transitively on $\Ga(\a)$, see \cite[Lemma 2.5]{567} for example. Thus $N_\a^{\Ga(\a)}$ is a transitive normal subgroup of $G_\a^{\Ga(\a)}$.
By Lemma \ref{imp-normal-sub}, $\soc(N_\a^{\Ga(\a)})=\soc(G_\a^{\Ga(\a)})$ and one of the following holds:
\begin{itemize}
        \item[(i)]  $N_\a^{\Ga(\a)}$ is a primitive permutation on $\Ga(\a)$;

        \item[(ii)] $N_\a^{\Ga(\a)}=\ZZ_p^k{:}\ZZ_m$ and $|\Ga(\a)|=p^k$, where $k\geqslant 2$, $m\div (p^d-1)$ for some divisor $d$ of $k$ with $d<k$;
        \item[(iii)] $N_\a^{\Ga(\a)}=\ZZ_3^4{:}\Q_8$ and $|\Ga(\a)|=3^4$.
\end{itemize}

\begin{lemma}\label{imp-stab}
Assume that  $G$ is $2$-arc-transitive on $\Ga$,  and $N\unlhd G$
with $N_\a\ne 1$ for $\a\in V$. Suppose that  $N_\a^{\Ga(\a)}$ is not primitive on $\Ga(\a)$. Then one of the following holds:
 \begin{itemize}
        \item[(1)] $N_\a=\ZZ_p^k{:}(\ZZ_{m_1}{.}\ZZ_m)=(\ZZ_p^k\times \ZZ_{m_1}){.}\ZZ_m$, $|\Ga(\a)|=p^k$ and $N_\a^{[1]}\cong \ZZ_{m_1}$, where $m_1\div m$, $m\div (p^d-1)$ for some divisor $d$ of $k$ with $d<k$;
        \item[(2)] $N_\a=\ZZ_3^4{:}(Q.\Q_8)=(\ZZ_3^4\times Q){.}\Q_8$, $|\Ga(\a)|=3^4$ and $Q\cong N_{\a}^{[1]}$, where $Q$ is isomorphic to a subgroup of $\Q_8$.
\end{itemize}
\end{lemma}
\proof
By the foregoing argument, we may let $N_\a^{\Ga(\a)}=KH$, where
$K=\soc(N_\a^{\Ga(\a)})\cong \ZZ_p^k$, and either $H\cong \ZZ_m$ or $p^k=3^4$ and $H\cong \Q_8$.
Without loss of generality, let $H=(N_\a^{\Ga(\a)})_\b$ for some $\b\in \Ga(\a)$.
Then $N_\a^{[1]}/N_{\a\b}^{[1]}$ is isomorphic a normal subgroup of $H$, see (\ref{ext-2}) given in Section \ref{sect=local}.

\vskip 5pt

Assume first that $p^k=4$. In this case, we have $H=1$ and $N_\a^{\Ga(\a)}\cong \ZZ_2^2$, and so  $N_\a$ acts faithfully on $\Ga(\a)$, refer to \cite[Lemma 2.3]{567}. Then $N_\a= \ZZ_2^2$, desired as part (1) of this lemma.


\vskip 5pt

Now assume that $p^k\ne 4$. Then $|\Ga(\a)|=p^k>5$.  By \cite[Theorem 4.7]{Weiss}, $G_{\a\b}^{[1]}=1$, and so $N_{\a\b}^{[1]}=1$, where $\b\in \Ga(\a)$.
Then $N_\a^{[1]}$ is isomorphic a normal subgroup of $H$. It is easily shown that
$|\Aut(N_\a^{[1]})|<p^k$.
Let $P$ be a Sylow $p$-subgroup of $N_\a$. Then
$PN_\a^{[1]}/N_\a^{[1]}$ is the unique Sylow $p$-subgroup of $N_\a/N_\a^{[1]}$, in particular, $PN_\a^{[1]}\unlhd N_\a$. Noting that $PN_\a^{[1]}/\C_{PN_\a^{[1]}}(N_\a^{[1]})$ is isomorphic to a subgroup of $\Aut(N_\a^{[1]})$, it follows that $p$ is a divisor of $|\C_{PN_\a^{[1]}}(N_\a^{[1]})|$. Let $Q$ be a Sylow $p$-subgroup of $\C_{PN_\a^{[1]}}(N_\a^{[1]})$. Then $Q$ is characteristic in $\C_{PN_\a^{[1]}}(N_\a^{[1]})$, and hence $Q$ is normal in $N_\a$. This implies that
$\O_p(N_\a)\ne 1$, where $\O_p(N_\a)$ is the maximal normal $p$-subgroup of $N_\a$. Since $N_\a\unlhd G_\a$, we have $\O_p(N_\a)\unlhd G_\a$.
Then either $\O_p(N_\a)\leqslant G_\a^{[1]}$ or $\O_p(N_\a)$ acts transitively on $\Ga(\a)$.
Noting that $(p,|N_\a^{[1]}|)=1$, we know that  $\O_p(N_\a)$ is faithful and transitive on $\Ga(\a)$. It follows that $P=\O_p(N_\a)\cong \ZZ_p^k$, and thus $N_{\a}=P{:}N_{\a\b}=(P\times N_\a^{[1]}){.}H$. Then part (1) or (2) of the lemma follows.
\qed

\vskip 10pt

The following theorem together with Theorem \ref{basic-gph-1} and Lemma \ref{imp-stab}   fulfil  the proof of Theorem \ref{basic-gph}.

\begin{theorem}\label{PA-type}
Assume that $\Ga=(V,E)$ is a basic $2$-arc-transitive graph with respect to  $G$, and  $G^*=\l G_\a,G_\b\r$ for and edge $\{\a,\b\}\in E$.
If $\Ga$ is not a complete bipartite graph then either  $\soc(G^*)$ is a nonabelian simple group, or every simple direct factor of $\soc(G^*)$ is semiregular on $V$.
\end{theorem}
\proof
Assume that $\Ga$ is not a complete bipartite graph.
Then $G^*$ is faithful on each of its orbits on $V$.
In view of Theorem \ref{basic-gph-1}, we may assume that
$G^*$ is a quasiprimitive permutation group of type {\rm PA} on each $G^*$-orbit   on $V$.
Let $N=\soc(G^*)$, and write
$N=T_1\times\cdots\times T_l$, where $l\geqslant 2$ and $T_i$ are isomorphic nonabelian simple groups. Then $N_\a\ne 1$, and $N_\a$ has no composition factor isomorphic to $T_1$, see
\cite[III(b)(i)]{Praeger-ON}.

\vskip 5pt

Let $U$ be a $G^*$-orbit on $V$ with $\a\in U$, and let $W=V\setminus U$ if $\Ga$ is bipartite. Clearly, $U$ is an $N$-orbit, and if $\Ga$ is bipartite then $W$ is also an $N$-orbit.
Note that $N$ is a minimal normal subgroup of both $G$ and $G^*$.
Since $G^*=NG_\gamma$ for $\gamma\in V$, it follows that both
$G$ and $G_\gamma$ act transitively on $\Omega:=\{T_1,\ldots,T_l\}$ by conjugation.
Let \[\CC_\gamma=\{(T_i)_\gamma\mid 1\leqslant i\leqslant l\},\,\,
\CC=\cup_{\gamma\in V}\CC_\gamma.\]
For $1\leqslant i\leqslant l$ and $x\in G$, we have $T_i^x\in \Omega$, and so
\[(T_i)_\gamma^x=(T_i\cap G_\gamma)^x=T_i^x\cap G_{\gamma^x}=(T_i^x)_{\gamma^x}\in \CC, \, \forall \gamma\in V.\]
It follows that  $G_\gamma$ acts transitively on $\CC_\gamma$  by conjugation, and
$\CC$ is a conjugacy class of subgroups in $G$.
In particular,
all orbits of each $T_i$ on $V$ have the same length $|T_1:(T_1)_\a|$.

\vskip 5pt

{\bf Case 1}.
Assume that  $N_\a^{\Ga(\a)}$ is primitive on $\Ga(\a)$.
For any $\gamma\in V$, letting $\gamma=\a^g$ for some $g\in G$, we have
\[\Ga(\gamma)=\Ga(\a)^g,\, N_\gamma=N\cap G_{\a^g}=(N\cap G_\a)^g=N_\a^g.\]
It follows that $N_\gamma$ acts primitively on $\Ga(\gamma)$. Thus $N$ is locally-primitive on $\Ga$.
If $T_1$ is intransitive on every $G^*$-orbit, then $T_1$ is
semiregular on $V$, see \cite[Lemma 2.6]{567}, and the result is true.
Suppose that $T_1$ is transitive on one of the $G^*$-orbits, say $U$. Since $T_l$ centralizes $T_1$, by \cite[Theorem 4.2A]{Dixon},  $T_l$ is semiregular on $U$.
This implies that both $T_1$ and $T_l$ are regular on $U$. Then $N=T_lN_\a$, and
so \[T_1\times\cdots\times T_{l-1}\cong N/T_l=T_lN_\a/T_l\cong N_\a/(T_l\cap N_\a)=N_\a/(T_l)_\a.\] Thus
$N_\a$ has a composition factor isomorphic to $T_1$, a contradiction.

\vskip 5pt

{\bf Case 2}.
Assume that $(T_1)_\a\leqslant G_\a^{[1]}$. Then $(T_1)_\a\leqslant (T_1)_\b$, where $\b\in \Ga(\a)$.
Recalling that $\CC$ is a conjugacy class of subgroups in $G$, it follows that $|(T_1)_\gamma|=|(T_1)_\a|$  for any $\gamma\in V$. Thus we have $(T_1)_\a= (T_1)_\b$. Since $(T_1)_\b\unlhd N_\b$ and $N_\b$ acts transitively on $\Ga(\b)$,
all $(T_1)_\a$-orbits on $\Ga(\b)$ have the same length. It follows that
$(T_1)_\a$ fixes $\Ga(\b)$ point-wise, i.e., $(T_1)_\b=(T_1)_\a\leqslant G_\b^{[1]}$. It follows from the connectedness of $\Ga$  that $(T_1)_\gamma=(T_1)_\a$  for any $\gamma\in V$. This forces that $(T_1)_\a=1$, and then our result is true in this case.

\vskip 5pt

{\bf Case 3}.
Assume that  $N_\a^{\Ga(\a)}$ is not primitive on $\Ga(\a)$,
and $(T_1)_\a\not\leqslant G_\a^{[1]}$.
Recall that  $G_\a$ acts transitively on $\CC_\a$. This implies that
$G_\a$ acts transitively on $\{(T_1)_\a^{[1]},\ldots,(T_l)_\a^{[1]}\}$,
$(T_1)_\a\times \cdots\times (T_l)_\a\unlhd G_\a$, and
$(T_i)_\a\not\leqslant G_\a^{[1]}$ for $1\leqslant i\leqslant l$.
By Lemma \ref{imp-normal-sub}, we have that
\[\soc(((T_1)_\a\times \cdots\times (T_l)_\a)^{\Ga(\a)})=\soc(G_\a^{\Ga(\a)})=\soc(N_\a^{\Ga(\a)})\cong \ZZ_p^k,\] and a Sylow $p$-subgroup of $N_\a$ has order $p^k$, where $p$ is a prime $p$ and $k\geqslant 2$.
By Lemma \ref{imp-stab}, $N_\a^{[1]}$ has order coprime to $p$, and
thus $(p,(T_i)_\a^{[1]})=1$ for $1\leqslant i\leqslant l$.

\vskip 5pt

Let $P_i$ be a Sylow $p$-subgroup of $(T_i)_\a$, where $1\leqslant i\leqslant l$. Then $P=P_1\times\cdots\times P_l$ is a Sylow $p$-subgroup of $N_\a$, and thus \[P\cong P^{\Ga(\a)}=\soc(N_\a^{\Ga(\a)})=\soc(G_\a^{\Ga(\a)}),\] and $\O_p((T_i)_\a^{\Ga(\a)})=P_i^{\Ga(\a)}\cong P_i$ for each $i$.
It follows that \[\soc(G_\a^{\Ga(\a)})=P_1^{\Ga(\a)}\times \cdots\times P_l^{\Ga(\a)}.\]
Let $K_i$ be the preimage of $P_i^{\Ga(\a)}$ in $(T_1)_\a\times \cdots\times (T_l)_\a$.
Then $K_i=(T_i)_\a^{[1]}P_i$ for $1\leqslant i\leqslant l$. It is easily shown
that $G_\a$ acts transitively on $\{K_1,\ldots,K_l\}$ by conjugation.
Then $G_\a^{\Ga(\a)}$ acts transitively on $\{P_1^{\Ga(\a)},\ldots,P_l^{\Ga(\a)}\}$ by conjugation, which is impossible by Lemma \ref{affine-2trans}.
This completes the proof of the theorem.
\qed

\vskip 20pt

Finally, we give a proof of Theorem {\rm \ref{bi-prim}}.

\vskip 5pt

\noindent {\it Proof of Theorem {\rm \ref{bi-prim}}}.
Assume that $G$ is a $2$-arc-transitive group of $\Ga=(V,E)$. Let $G^*=\l G_\a,G_\b\r$ for  $\{\a,\b\}\in E$. If $\Ga$ is not bipartite and
$G$ is primitive on $V$ then $\soc(G)$ is either simple or regular on $V$ by
\cite[Theorem A]{Prag-prim-2arc}, and the result is true.

\vskip 5pt

Assume next that $\Ga$ is a bipartite graph with two parts $U$ and $W$,
and that $G^*$ acts primitively on both $U$ and $W$. If $G^*$ is unfaithful on
$U$ or $W$ then $\Ga$ is a complete bipartite graph. Thus we assume further that
$G^*$ is faithful on both $U$ and $W$. Let $\a\in U$ and $\b\in W$.

\vskip 5pt

{\bf Case 1}. Assume that $\soc(G)\leqslant G^*$. Then $\Ga$ is a basic $2$-arc-transitive graph with respect to $G$. By Theorem \ref{basic-gph-1}, $\soc(G)=\soc(G^*)$, and either part (2) of Theorem \ref{bi-prim} holds or $G^*$ is a primitive permutation group of type {\rm PA} on $U$. For the latter case,
every simple direct factor of $\soc(G^*)$ is not semiregular on $U$, refer to \cite[page 391, III(b)(i)]{Lie-P-S}. Then part (2) of Theorem \ref{bi-prim} occurs
by  Theorem \ref{PA-type}.

\vskip 5pt

{\bf Case 2}. Assume that $\soc(G)\not\leqslant G^*$. Let $M$ be a minimal normal subgroup of $G$ with $M\not\leqslant G^*$. Then, noting that $|G:G^*|=2$, we have
$G=G^*\times M$ and $|M|=2$. This implies that $\soc(G)=\soc(G^*)\times M$. Set $M=\l x\r$. Then $G_{\a^x}=G_\a^x=G_\a$, and
so $G_\a$ acts $2$-transitively on $\Ga(\a^x)$. Considering the (faithful) action of $G^*$ on $U$, by \cite[Theorem A]{Prag-prim-2arc}, $\soc(G)^*$ is either simple or regular on $U$. Similarly, $\soc(G)^*$ is either simple or regular on $W$.
Then part (3) of Theorem \ref{bi-prim} follows.
 \qed

\end{document}